# Estrategias para cubrir la demanda insatisfecha en una fábrica alimenticia

Dr. Carlos Eduardo Mendoza Durán[1], Francisco Gutiérrez Albitos[2],
Rodrigo Leo Paniagua[3], Héctor Manuel Rebollo Arana[4].

*Resumen* — **La insatisfacción de la demanda puede ser un problema grave en una empresa dedicada a la producción de bienes de consumo, y la atribución de dicha insatisfacción a una causa específica es una etapa crucial en la solución del problema. En este trabajo se presenta un caso de estudio cuya solución fue encontrada mediante la formulación de un modelo de programación no lineal.**

*Palabras clave* — **demanda insatisfecha, investigación de operaciones, programación matemática, modelos de pronóstico, control de piso de producción**

## Introducción

Este trabajo tiene como propósito reportar la metodología empleada en un proyecto terminal de licenciatura en ingeniería industrial, cuyo objetivo fue determinar las causas por las cuales una planta dedicada a la producción de pasta alimenticia de trigo no era capaz de suministrar la totalidad de sus pedidos.

*Antecedentes*

La planta donde se desarrolló el trabajo se encuentra ubicada en la Ciudad de México, y produce alrededor de 2,500 toneladas de pasta de sémola de trigo al mes. Sus clientes principales son tiendas de autoservicio, para las cuales produce paquetes de pasta que se venden bajo marcas libres, aunque cuenta también con dos marcas propias. Al momento de realizar la visita de diagnóstico a la empresa, los directivos reportaron que según sus estimaciones, aproximadamente el 10% de los pedidos que ordenan los clientes no se entregan, presumiblemente debido a que el espacio disponible para almacenar el producto terminado es muy pequeño, y no por un déficit de capacidad de producción en la planta.

Con la finalidad de que la empresa sea capaz de cubrir mejor su demanda, se plantearon los siguientes objetivos:
1. Cuantificar el tamaño de la demanda que no se satisface.
2. Determinar la viabilidad de pronosticar la demanda futura mediante modelos sencillos de series de tiempo.
3. Establecer la causa principal de la insatisfacción de la demanda evaluando lo siguiente:
    a. Si la capacidad productiva real de la planta es suficiente para suministrar todos los pedidos.
    b. Si una mejora en el control de piso de producción puede lograr la satisfacción de la demanda.
    c. Si el tamaño del almacén de producto terminado actual es una restricción significativa en la satisfacción de la demanda.

*Metodología*

La metodología del proyecto fue la siguiente:
1. Recolección de información histórica acerca de la demanda y las ventas.
2. Realización de un análisis comparativo entre pedidos recibidos y ventas realizadas para estimar la demanda no satisfecha.
3. Elaboración de modelos de pronóstico de series de tiempo (media, media móvil, suavizado exponencial, regresión lineal simple) para la demanda, y evaluación de su fiabilidad.
4. Medición de la capacidad productiva de la planta, y determinación del impacto de los tiempos de arranque de las máquinas en dicha capacidad.
5. Elaboración de un modelo matemático de control de piso de producción que considere las restricciones debidas tanto a la capacidad de producción como a la capacidad de almacenamiento de productos terminados.

---

[1] El Dr. Carlos Mendoza es coordinador académico de la maestría en logística y del doctorado en ingeniería industrial de la Universidad Anáhuac México Norte, así como profesor de planta en la Facultad de Ingeniería y en el Centro de Alta Dirección en Ingeniería y Tecnología de la misma universidad.
[2] Francisco Gutiérrez es estudiante de noveno semestre de ingeniería industrial en la Universidad Anáhuac México Norte.
[3] Rodrigo Leo es estudiante de noveno semestre de ingeniería industrial en la Universidad Anáhuac México Norte.
[4] Héctor Rebollo es estudiante de noveno semestre de ingeniería industrial en la Universidad Anáhuac México Norte.

**Análisis de la demanda**

*Determinación del tamaño de la demanda no satisfecha*

En primer lugar se recolectó la información de las ventas y pedidos realizados cada mes entre enero de 2011 y diciembre de 2013. Dicha información se presenta como una serie temporal en la figura (1), donde se observa que en la mayoría de los meses las ventas son inferiores a la demanda, lo cual confirma el hecho de que hay demanda insatisfecha. Se estimó el tamaño de dicha demanda insatisfecha ($D_{\text{ins}}$) como el promedio de los porcentajes de pedidos no entregados cada mes entre los años 2011 y 2013, resultando en un 13.34%.

$$D_{\text{ins}} = \frac{1}{36} \sum_{i=1}^{36} \frac{D_i - V_i}{D_i} = 13.34\%$$

Donde $D_i$ y $V_i$ son la demanda en kilogramos y la venta en kilogramos (respectivamente) del periodo (mes) $i$.

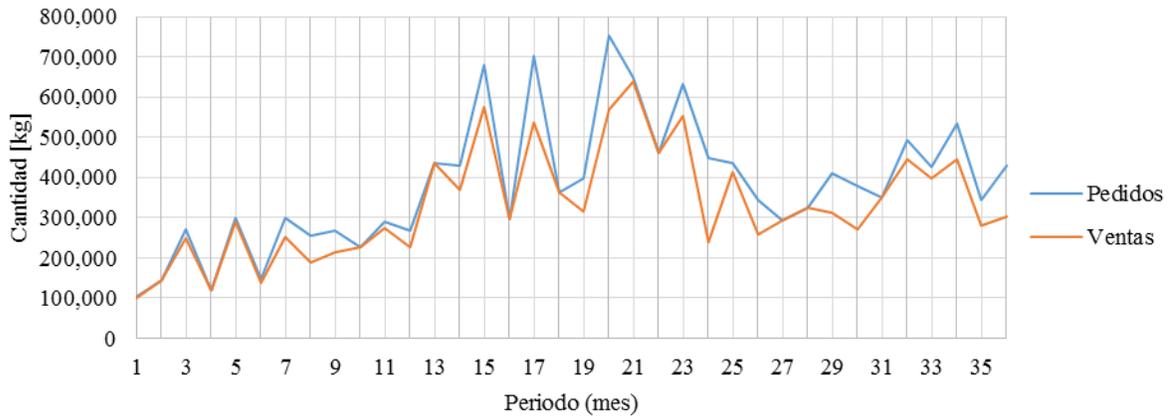

**Figura 1.** Pedidos y ventas (en kilogramos) del periodo 2011-2013.

En la figura (1) Se puede observar que existe demanda insatisfecha independientemente de la cantidad de kilogramos que se piden en el mes[5], por lo que es razonable suponer que la razón por la cual los pedidos no se satisfacen no es falta de capacidad productiva en la planta. Sin embargo, esto se analizará con detalle en la siguiente sección.

*Desarrollo de modelos de pronóstico para la demanda*

Con el fin de pronosticar la demanda, se evaluaron distintos modelos de pronóstico basados en series de tiempo: media simple, medias móviles (variando el número de periodos entre 2 y 23), modelos de suavizado exponencial (simple, doble –con tendencia– y triple –con tendencia y estacionalidad–) y regresión lineal. En el cuadro (1) se muestran los modelos con el MAPE[6] correspondiente.

La demanda presenta un comportamiento normal con media $\mu = 380{,}506$ y desviación estándar $\sigma = 160{,}689$ con un nivel de significancia $\alpha = 0.05$ contra un valor $p = 0.121$ correspondiente a la prueba de Anderson-Darling. El coeficiente de variación de la demanda es de 42.23%. Debido a esta razón y a la notable alternancia de crestas y valles sin periodicidad notable, los pronósticos mediante series de tiempo no pueden considerarse suficientes para predecir los pedidos futuros (Ramesh, 2009). Sin embargo, bajo la hipótesis de normalidad –ya probada–, pueden construirse estimadores estadísticamente fiables para los tamaños de futuros pedidos.

---

[5] Por ejemplo, en el mes 7 se observa una insatisfacción de aproximadamente 50,000 kg, mientras que la demanda del periodo es de 300,000 kg. Por otra parte, el mes 23 la demanda es de más de 600,000 kg y la insatisfacción es igualmente de 50,000 kg.

[6] El MAPE (*mean average percentage error*) es el promedio de los errores porcentuales absolutos y una buena medida de error para comparar distintos modelos de pronóstico (Krakewski et al. 2008).

| Modelo | MAPE |
|---|---|
| Promedio simple | 28.3% |
| Promedio móvil de 17 periodos[7] | 17.3% |
| Suavizado exponencial[8] simple con $\alpha = 0.22$ | 25.0% |
| Suavizado exponencial doble (modelo de Holt) con $\alpha = 0.42$ y $\gamma = 0.07$ | 30.9% |
| Suavizado exponencial triple (modelo de Winters) con $\alpha = 0.42$, $\gamma = 0.07$ y $\delta = 0.2$ | 27.5% |
| Regresión lineal | 32.0% |

**Cuadro 1.** Modelos de pronóstico y MAPEs correspondientes.

### Análisis de capacidad

El objetivo del análisis de capacidad es determinar si la capacidad productiva de la planta (medida en kg/mes) es suficiente para cumplir con las órdenes de manufactura que se reciben de parte de los clientes. Si bien en la figura (1) se observa que la insatisfacción de pedidos no parece depender de la demanda del periodo en cuestión, es necesario realizar un análisis específico encaminado a determinar esto con mayor certeza.

La planta cuenta con cuatro máquinas, cada una con diferentes tasas nominales de producción $t_j$ (medidas en kg/h), así como un índice de eficiencia $e_j$ y un porcentaje de merma $m_j$ debida al proceso de extrusión de la pasta cruda. La tasa neta $\tau_j$ se calcula entonces como $\tau_j = t_j e_j (1 - m_j)$. Al momento de comenzar a procesar un pedido en una máquina en particular, la máquina en cuestión requiere un tiempo de arranque que oscila alrededor de doce horas. Este tiempo de arranque considera el cambio de molde de extrusión, la limpieza de la máquina y el establecimiento de las condiciones de calor y humedad necesarias para la cocción y el secado de la pasta.

Para estimar el impacto debido a $n$ arranques es necesario considerar un estimador razonable de la pérdida de capacidad ($\Delta C(n)$): el promedio de las tasas reales $\tau_j$ multiplicado por el tiempo de arranque promedio $s_j$ multiplicado por el número de arranques $n$:

$$\Delta C(n) = n \frac{1}{4} \left( \sum_{j=1}^{4} \tau_j \times \sum_{j=1}^{4} s_j \right)$$

De esta forma se puede calcular la capacidad neta de la planta como su capacidad nominal mensual $C_N = 176 \sum_{j=1}^{4} \tau_j$ menos el cambio en capacidad $\Delta C(n)$ debido a $n$ arranques. En la figura (2) se muestra un gráfico histórico de la capacidad utilizada (ventas) superpuesta a la capacidad requerida (demanda). Excepto en un periodo, la capacidad utilizada es menor o igual a la capacidad requerida. Sin embargo, para completar el análisis, se debe considerar también la capacidad *disponible*: $C_N - \Delta C(n)$.

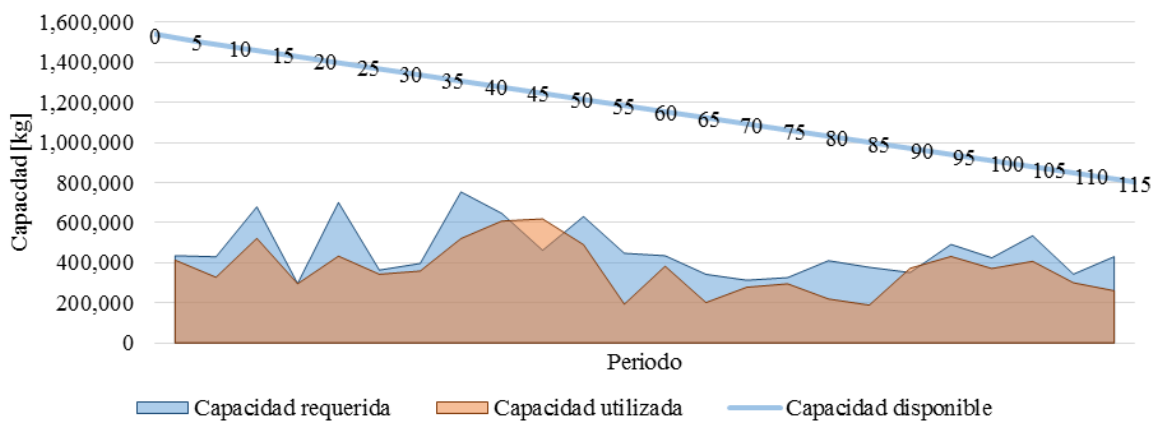

**Figura 2.** Gráfico de capacidad utilizada, requerida y disponible (en kilogramos).

Como se puede observar, la capacidad disponible depende linealmente del número de arranques efectuados en la planta, y se muestra en la figura (2) como una línea recta azul, cuya altura corresponde a la capacidad en función del número de arranques (mostrado directamente sobre la línea). Se observa que inclusive con cien arranques mensuales, la capacidad disponible excede la capacidad requerida. Y cien arranques mensuales equivalen en promedio a cinco arranques diarios, lo cual es excesivo considerando que hay cuatro máquinas, que cada arranque tarda alrededor de doce horas y que la jornada laboral en la planta es de ocho horas.

---

[7] De todos los modelos de media móvil, se muestra sólo el caso de 17 periodos porque es el que mejor ajuste muestra basado en el MAPE.

[8] Para las series basadas en suavizado exponencial, las constantes de suavización se eligieron como aquellas que minimizaron el MAPE.

Por lo tanto, debe concluirse que la insatisfacción de la demanda no es directamente atribuible a una falta de capacidad productiva en la planta.

### Capacidad de almacenamiento y control de la producción

En esta sección se desarrolla un modelo de programación matemática que pretende determinar si es factible cubrir completamente la demanda de la planta con las restricciones actuales de capacidad de producción e inventario. Dicho modelo se desarrolla de forma general para posteriormente aplicarlo al problema de la planta. Es conveniente mencionar que originalmente el modelo fue llevado a cabo en dos etapas. En la primera etapa, se consideró como única restricción de producción la capacidad productiva de la planta, y los resultados mostraron que la demanda podía satisfacerse simplemente distribuyendo correctamente los pedidos en las máquinas disponibles. En la segunda etapa se agregó como restricción la capacidad de almacenamiento de los inventarios de producto terminado. Los resultados de la solución de este modelo extendido se muestran y analizan al final de esta sección.

*Enunciado del problema*

En la planta se tienen, en cierto momento, $n$ pedidos que se traducirán en órdenes de manufactura. Cada pedido $i$ tiene asociados los siguientes parámetros: un tamaño de pedido (en kilogramos) denotado $Q_i$ y una fecha de entrega límite (en días) denotada $E_i$. En la fábrica se cuenta con $m$ máquinas. Cada máquina $j$ tiene asociada una tasa de producción (en kilogramos sobre hora) denotada $\tau_j$ y un tiempo de arranque medio (en horas) denotado $s_j$. Finalmente, la fábrica cuenta con un almacén de producto terminado cuya capacidad máxima (en kilogramos) es $A$. Se considera que cada día laboral tiene $h$ horas netas de trabajo disponibles. El problema consiste en encontrar cómo deben distribuirse las órdenes en cada una de las máquinas de forma que: se produzca la mayor cantidad de pasta posible[9], todas las órdenes sean entregadas a tiempo (es decir, tengan una holgura[10] positiva), y la capacidad máxima del almacén nunca se sobrepase.

*Construcción del modelo*

En primer lugar se definen $nm$ variables binarias que servirán como las variables de decisión del problema:
$$x_{i,j} = \begin{cases} 1 & \text{si la orden } i \text{ se produce en la máquina } j \\ 0 & \text{en otro caso} \end{cases}$$

Debido a que cada orden se procesa a lo más en una máquina, cada vector $(x_{i,1}, x_{i,2}, \ldots, x_{i,m})$ deberá tener a lo más una entrada igual a uno, y las demás deberán ser iguales a cero. Esto constituye el primer conjunto de restricciones para el problema:

$$\forall i \in \{1, \ldots, n\}, \quad \sum_{j=1}^{m} x_{i,j} \leq 1 \tag{1}$$

A continuación se define una matriz $T$ de dimensión $n \times m$ –llamada *matriz de tiempos de proceso*– que contiene en la entrada $i, j$ el tiempo de proceso para el pedido $i$ en la máquina $j$:

$$T = \begin{pmatrix} x_{1,1}\left(s_1 + \dfrac{Q_1}{\tau_1}\right) & \cdots & x_{1,m}\left(s_m + \dfrac{Q_1}{\tau_m}\right) \\ \vdots & \ddots & \vdots \\ x_{n,1}\left(s_1 + \dfrac{Q_n}{\tau_1}\right) & \cdots & x_{n,m}\left(s_m + \dfrac{Q_n}{\tau_m}\right) \end{pmatrix}$$

La suma de los valores de la fila $i$ de la matriz $T$ es igual al tiempo de proceso requerido para fabricar la orden $i$. Este tiempo de proceso se denota $T_{P_i}$:

$$T_{P_i} = \sum_{j=1}^{m} [T]_{i,j} = \sum_{j=1}^{m}\left[x_{i,j}\left(s_j + \dfrac{Q_i}{\tau_j}\right)\right] \tag{2}$$

---

[9] Se asume, con fundamento en el análisis de la demanda, que la producción de esta planta es igual a su venta. Por esta razón, la cantidad de pasta producida (función objetivo) puede considerarse proporcional a la utilidad, descartando las variaciones de precio entre productos.

[10] La holgura de una orden se define como la diferencia entre su fecha límite de entrega y su fecha de terminación. En otras palabras, se trata de la cantidad de tiempo que la orden ya finalizada «espera» antes de ser entregada.

El siguiente paso en la construcción del modelo es determinar el momento en que la orden $i$ será terminada. Este tiempo requiere dos funciones auxiliares: el tiempo de inicio de la orden $i$ en la máquina $j$, denotado $I_{i,j}$; y el tiempo de finalización de la orden $i$ en la máquina $j$, denotado $F_{i,j}$. Debido a razones de espacio, se omite aquí el desarrollo detallado de dichas funciones. En general, el tiempo de finalización $F_{i,j}$ se escribe como sigue (considerando que para cualquier $j$ se debe cumplir que $I_{i,j} = F_{(i-1),j}$):

$$F_{i,j} = \begin{cases} \sum_{j=1}^{m}\left[x_{1,j}\left(s_j + \frac{Q_1}{\tau_j}\right)\right] & \text{si } i = 1 \\ I_{i,j} + \sum_{j=1}^{m}\left[x_{i,j}\left(s_j + \frac{Q_i}{\tau_j}\right)\right] & \text{en otro caso} \end{cases} \quad (3)$$

A continuación se escribe la fecha de finalización de la orden $i$, denotada $L_i$, como el tiempo de finalización de la orden en la máquina que le corresponde[11]:

$$L_i = \sum_{j=1}^{m} x_{i,j} F_{i,j} \quad (4)$$

Con la fecha de finalización puede formularse la restricción principal del problema: la diferencia entre la fecha de entrega $E_i$ y la fecha de finalización $L_i$ debe ser mayor que cero para que la orden se entregue a tiempo. Esta diferencia se conoce como *holgura de orden* y se denota $H_i$. La restricción se escribe como:

$$\forall i \in \{1, \ldots, n\}, \quad H_i = hE_i - L_i > 0 \quad (5)$$

Para plantear las restricciones debidas a la capacidad del almacén se define la *interacción* de dos órdenes de la forma siguiente: se dice que la orden $i$ interactúa con la orden $k$ si se cumple alguna de las dos siguientes condiciones: (1) la orden $k$ entró al almacén después de la entrada de la orden $i$ y antes de la salida de la orden $i$ ($L_i < L_k < E_i$); (2) la orden $k$ entró al almacén antes de la entrada de la orden $i$, y salió después de la entrada de la orden $i$ ($L_k < L_i < E_k$). Finalmente se define un factor de interacción, denotado $\lambda_{i,k}$:

$$\lambda_{i,k} = \begin{cases} 1 & L_i < L_k < E_i \ \lor \ L_k < L_i < E_k \\ 0 & \text{en otro caso} \end{cases} \quad (6)$$

Mediante el factor de interacción descrito en la ecuación (6) se expresa la ocupación total del almacén (en kilogramos) mientras la orden $i$ está en él. Esta ocupación se denota $O_i$, y se restringe de la siguiente forma:

$$\forall i \in \{1, \ldots, n\}, \quad O_i = \left(\sum_{j=1}^{m} x_{i,j}\right)\left(\sum_{k=1}^{n}[\lambda_{i,k}Q_k] + Q_i\right) \leq A \quad (7)$$

Finalmente, se define la función objetivo $Z$ del problema como la cantidad (en kilogramos) de pasta producida:

$$Z = \sum_{i=1}^{n}\left[Q_i \sum_{j=1}^{m} x_{i,j}\right] \quad (8)$$

*Naturaleza del problema de programación*

Si bien la función objetivo presentada en la ecuación (8) es lineal, el modelo no lo es, ya que la ecuación (7) viola el supuesto de aditividad de la programación lineal (Hillier, 1991), por lo que no puede ser resuelto con un algoritmo PL (simplex o punto interior, por ejemplo). Por lo tanto, el problema pertenece al ámbito de la programación no lineal. Se puede verificar fácilmente que tampoco se trata de un problema de programación convexa, ya que la región definida por la restricción (1) es no convexa[12]. Por otra parte, está demostrado que el problema general de control de piso es un problema NP-completo para $m > 2$ (Garey, 1976), por lo que no es resoluble en tiempo polinómico[13]. Otros trabajos de investigación (Biegel, 1971) muestran que el tiempo de solución del problema de asignación en $m$ máquinas crece de forma proporcional a $m!$. La proposición de un algoritmo para la resolución de un problema de estas características queda fuera del alcance de este trabajo.

---

[11] En la expresión de la ecuación (4) se introduce el factor $x_{i,j}$ debido a que no se sabe de antemano en qué máquina se producirá cada orden.

[12] Esta afirmación se prueba comprobando que la matriz hessiana $Hf$ de $f(i) = \sum_{j=1}^{m} x_{i,j}$ es indefinida (Hillier, 1991) porque es la matriz nula: $[Hf]_{i,j} = 0$.

[13] A menos, por supuesto, que la hipótesis P=NP sea verdadera (véase Fortnow, 2009).

*Solución*

El modelo se implementó en Microsoft Excel® 2013 y se solucionó mediante el motor de gradiente reducido generalizado (GRG) del complemento Solver®. Debido a las características de este algoritmo, no puede asegurarse que los óptimos encontrados sean globales, pero al menos son puntos que satisfacen las condiciones de Karush-Kuhn-Tucker (Frontline, s.f), asegurando que la solución constituye un óptimo local.

*Análisis del aumento de la capacidad de almacenamiento de producto terminado*

Se contemplan dos escenarios para la capacidad de almacén: (a) una capacidad máxima de almacén de 84 pallets (capacidad actual) y (b) una capacidad máxima de almacén de 144 pallets (que puede alcanzarse instalando tres pisos de anaqueles en el almacén actual). El modelo se ejecutó para la demanda de cada mes del año 2013 (segmentada en pedidos), manteniendo la capacidad productiva de la planta constante y variando la capacidad de almacén. Los resultados del cuadro 2 muestran la clara reducción de la demanda insatisfecha.

| Mes | Número de pedidos no satisfechos $A = 84$ pallets | $A = 144$ pallets | Demanda [kg] | Mes | Número de pedidos no satisfechos $A = 84$ pallets | $A = 144$ pallets | Demanda [kg] |
|---|---|---|---|---|---|---|---|
| 1 | 1 | 0 | 435,536 | 7 | 0 | 0 | 351,338 |
| 2 | 0 | 0 | 342,621 | 8 | 2 | 0 | 491,975 |
| 3 | 0 | 0 | 294,082 | 9 | 1 | 0 | 424,908 |
| 4 | 0 | 0 | 326,342 | 10 | 2 | 0 | 535,150 |
| 5 | 1 | 0 | 410,814 | 11 | 0 | 0 | 343,411 |
| 6 | 1 | 0 | 377,721 | 12 | 1 | 0 | 430,795 |

**Cuadro 2.** Resultados de la solución encontrada para el modelo.

## Conclusión

La demanda promedio en 2013 fue de 397,058 kg, y se observa (cuadro 2) que los pedidos no pueden ser completamente satisfechos debido a limitantes de espacio a partir de una demanda de al menos 377,721 kg. Suponiendo que la demanda mensual sigue una distribución aproximadamente normal con media y desviación estándar estimadas a partir de la muestra ($\mu = 397,058$, $\sigma = 71,078$)[14], entonces la probabilidad de que la demanda sea mayor al límite crítico de 377,721 kg es de 60.72%. Por lo tanto, se puede esperar que más de la mitad de las veces la empresa no cumpla completamente con las entregas pactadas debido a limitaciones en el almacén actual de producto terminado.

Por otra parte, el valor crítico a partir del cual la demanda no se satisface por completo con una capacidad de almacén de producto terminado de 144 pallets es de 625,000 kg. La probabilidad de una demanda mensual mayor a dicho valor crítico es de 0.07%.

Por lo tanto, la ampliación de la capacidad del almacén de producto terminado de 84 pallets a 144 pallets mediante la instalación de anaqueles (que significa 70% más de espacio) redundará en un 137% más de satisfacción de la demanda, sin necesidad de instalar nuevas máquinas ni de buscar un nuevo local.

## Referencias

Biegel, J. *Operations planning, scheduling and control*. Industrial Engineering Handbook (editado por Maynard, H.). Tercera edición. McGraw-Hill. EUA (1971).

Fortnow, L. *The status of the P versus NP problem* (2009). Communications of the ACM (52), no. 9, pp. 78-86. doi: 10.1145/1562164.1562186.

Frontline Systems Inc. *Standard Excel Solver – GRG nonlinear solver stopping conditions*. Artículo de documentación de Solver disponible en http://www.solver.com/standard-excel-solver-grg-nonlinear-solver-stopping-conditions. Recuperado el 02/12/2014.

Garey, M., Johnson, D., Sethi, R. *The complexity of flowshop and jobshop scheduling* (1976). Mathematics of Operations Research 1 (2), pp. 117-119. doi: 10.1287/moor.1.2.117.

Hillier, F., Lieberman, G. *Introducción a la investigación de operaciones*. Quinta edición. McGraw-Hill. México (1991).

Krakewski L., Ritzman, L., Malhotra, M. *Administración de operaciones*. Octava edición. Pearson. México (2008).

Ramesh, V. *Effectively managing demand variability in CPG industry* (2009). Infosys White Papers.

---

[14] Esta suposición se confirmó mediante la prueba de normalidad de Anderson-Darling con un nivel de significancia $\alpha = 0.05$ contra un valor $p = 0.609$.